\documentclass[12pt]{article} 

\usepackage{amssymb}
\usepackage{amsmath}
\input{liemacs.sty}

\begin{document} 



\title{A non-smooth continuous unitary representation of a Banach--Lie group} 
\author{Daniel Belti\c t\u a\footnote{Institute of Mathematics ``Simion
Stoilow'' of the Romanian Academy, 
P.O. Box 1-764, RO-014700 Bucharest, Romania;  
email: \texttt{Daniel.Beltita@imar.ro}}\ \ 
and Karl-Hermann Neeb\footnote{Department of Mathematics, 
Darmstadt University of Technology, 
Schlossgartenstrasse 7, D-64289 Darmstadt, Germany;
email: \texttt{neeb@mathematik.tu-darmstadt.de}}}
\date{November 26, 2008}

\maketitle
\begin{abstract} In this note we show that the 
representation of the additive group of the Hilbert space $L^2([0,1],\R)$ 
on $L^2([0,1],\C)$ given by the multiplication operators 
$\pi(f) := e^{if}$ is continuous but its space of smooth 
vectors is trivial. This example shows that a continuous 
unitary representation of an infinite dimensional 
Lie group need not be smooth. 
\end{abstract} 


\section{Introduction} 

\begin{defn} Let $G$ be a Lie group modeled on a 
locally convex space (cf.\ \cite{Ne06} for a survey on locally convex 
Lie theory). 

Let $\cH$ be a complex Hilbert space and $U(\cH)$ be its 
unitary group. A  {\it unitary representation of $G$ on $\cH$} is a pair 
$(\pi, \cH)$, where $\pi \: G \to U(\cH)$ is a group homomorphism. 

A unitary representation $(\pi,\cH)$ is said to be {\it continuous} 
if the action $G \times \cH \to \cH, (g,v) \mapsto \pi(g)v$ is continuous. 
Since $G$ acts by isometries on $\cH$, it is easy to see that this 
condition is equivalent to the continuity of all orbit maps 
$\pi^v \: G \to \cH, g \mapsto \pi(g)v$. 

A unitary representation $(\pi,\cH)$ is said to be {\it smooth} 
if the space 
$$\cH^\infty := \{ v \in \cH \: \pi^v \in C^\infty(G,\cH) \} $$
of {\it smooth vectors} is dense. 
\end{defn}

Clearly, every smooth representation is continuous, and it is a natural 
question to which extent the converse also holds. 

\begin{rem} If $G$ is finite dimensional, then each continuous 
unitary representation is smooth. Even the subspace 
$\cH^\omega \subeq \cH^\infty$ of analytic vectors is dense 
(cf.\ \cite{Ga60}). 
\end{rem}

For the class of groups which are direct limits of 
finite dimensional Lie groups, 
Samoilenko's book \cite{Sa91} contains a variety of positive results 
on the existence of smooth vectors, in particular for abelian 
Lie groups, restricted direct products of $\SU_2(\C)$ and the group 
of infinite upper triangular matrices. 
More general existence 
results on differentiable vectors for limits of finite dimensional 
Lie groups can be found in \cite{Da96}.
See also \cite{Sh01} for existence of smooth vectors 
for particular classes of representations of diffeomorphism groups. 

However, 
the purpose of this note is to show that there is no 
automatic smoothness result for continuous unitary representations 
of infinite dimensional Lie groups. 
Even for the otherwise rather well-behaved class of abelian 
Banach--Lie groups. This will be shown by verifying that 
for the abelian Hilbert--Lie group $G = (L^2([0,1],\R),+)$, the 
unitary representation 
$$ \pi \: G \to U(L^2([0,1],\C)), \quad \pi(f)\xi := e^{if}\xi $$
is continuous, but its space $L^2([0,1],\C)^\infty$ of smooth vectors is trivial. 

Smoothness of a representation is a property that is crucial to make 
it accessible to Lie theoretic methods. In particular, for any smooth 
representation $(\pi, \cH)$ we obtain a representation of its 
Lie algebra $\g$ on the space $\cH^\infty$ of smooth vectors 
by skew-hermitian operators (cf.\ \cite{Ne08}). Our example 
shows that smoothness of a representation is an assumption that does 
not follow from continuity.

\section{The exponential representation}

\begin{prop} \label{prop:1} The unitary 
representation $(\pi,L^2([0,1],\C))$ of the additive group 
$G = L^2([0,1],\R)$, defined by $\pi(f)\xi = e^{if}\xi$, is continuous. 
\end{prop}

\begin{prf} First we observe that for any $t \in [0,1]$ and 
$f,g \in L^2([0,1],\R)$ we have 
$|e^{if(t)} - e^{ig(t)}| \leq |f(t)-g(t)|.$ 
For any $\xi \in L^2([0,1],\C) \cap L^\infty([0,1],\C)$ we thus obtain 
\begin{align*}
\|\pi(f)\xi - \pi(g)\xi\|_2^2 
&= \int_0^1 |e^{if(t)} - e^{ig(t)}|^2\cdot |\xi(t)|^2\, dt 
\leq \|\xi\|_\infty^2 \int_0^1 |f(t) - g(t)|^2\, dt \\
&=  \|\xi\|_\infty^2 \|f-g\|_2^2. 
\end{align*}
This implies that the orbit map $\pi(\cdot)\xi$ is continuous if 
$\xi$ is bounded, and since the set of bounded elements is dense 
in $L^2([0,1],\C)$, the continuity of $\pi$ follows. 
\end{prf}

To show that the space $L^2([0,1],\C)^\infty$ of smooth vectors is trivial, 
we put $\cH := L^2([0,1],\C)$ and consider the functions 
$$ f_\lambda(t) := |t-\lambda|^{-\frac{1}{4}}, \quad \lambda \in [0,1], $$
in $L^2([0,1],\R)$. For each $\lambda$, the continuous unitary representation 
$\pi$ defines a continuous unitary one-parameter group 
$$ \pi_\lambda(t)\xi := e^{it f_\lambda}\xi, $$
whose infinitesimal generator is the multiplication operator 
$$ M_\lambda \: \cD_\lambda \to \cH, \quad M_\lambda \xi := f_\lambda \xi, 
\quad \cD_\lambda := \{ \xi \in \cH \: \|f_\lambda \xi\|_2 < \infty\}. $$
In particular, the set of smooth vectors for this one-parameter group 
is the dense subspace 
$$ \cD_\lambda^\infty := \{ \xi \in \cH \: (\forall n \in \N) \ 
\|f_\lambda^n \xi\|_2 < \infty\}. $$
Therefore it remains to show that 
$\bigcap_{\lambda \in [0,1]} \cD_\lambda^\infty = \{0\}$. 

\begin{prop} \label{prop:2} If $\xi \in L^2([0,1],\C)$ has the property that 
$f_\lambda^4 \xi \in L^2([0,1],\C)$ holds for each 
$\lambda \in [0,1]$, then $\xi = 0$. 
\end{prop}

\begin{prf} Replacing $\xi$ by $|\xi|$, we may w.l.o.g.\ assume that 
$\xi \geq 0$. 

For $n \in \N$, let 
$M_n := \{ t \in [0,1] \: \frac{1}{n}\leq \xi(t) \leq n \}$ 
and note that $\xi = \lim_{n \to \infty} \xi \chi_{M_n}$ 
holds in $L^2([0,1],\C)$. If $f_\lambda^k \xi \in L^2([0,1],\C)$, then 
$f_\lambda^k \xi \chi_{M_n}\in L^2([0,1],\C)$ for any $k,n \in \N$ 
and hence  $f_\lambda^k \chi_{M_n} \in L^2([0,1],\C)$. 
We may therefore assume that $\xi = \chi_M$ is the characteristic 
function of some measurable subset $M \subeq [0,1]$. 

Suppose that $M$ has positive measure and 
that $f_\lambda^4 \xi \in L^2([0,1],\C)$ holds for each $\lambda \in [0,1]$. 
We have to show that this assumption leads to a contradiction. 
Let $\lambda \in M \cap ]0,1[$ be a Lebesgue point of $\xi$ 
(\cite[Thm.~7.11]{Ru86}), so that 
$$ 1 
= \lim_{h \to 0} \frac{|M \cap [\lambda-h,\lambda+h]|}{2h}  
= \lim_{h \to 0} \frac{1}{2h} \int_{\lambda - h}^{\lambda+h} \chi_M(t)\, dt. $$
We then find the two estimates 
$$ \int_{\lambda - h}^{\lambda+h} f_\lambda^4 \chi_M(t)\, dt 
\leq \|f_\lambda^4 \chi_M\|_2 \cdot \sqrt{2h} \to 0 
\quad \mbox{ as } \quad h \to 0, $$
and, likewise, for $h \to 0$, 
$$ \int_{\lambda - h}^{\lambda+h} f_\lambda^4 \chi_M(t)\, dt 
= \int_{\lambda - h}^{\lambda+h} \frac{1}{|t-\lambda|} \chi_M(t)\, dt 
\geq \frac{1}{h} \int_{\lambda - h}^{\lambda+h} \chi_M(t)\, dt 
\to 2. $$
These two estimates are contradictory, which completes 
the proof. 
\end{prf}

\begin{thm} The unitary representation of $G = L^2([0,1],\R)$ 
on $\cH = L^2([0,1],\C)$ defined by $\pi(f)\xi= e^{if}\xi$ 
is continuous, but all its smooth vectors are trivial.
\end{thm} 

\begin{prf} The continuity has been verified in Proposition~\ref{prop:1}. 
If $\xi \in \cH^\infty$ is a smooth vector, then 
$f_\lambda^4 \xi \in L^2([0,1],\C)$ holds for each $\lambda \in [0,1]$, 
so that Proposition~\ref{prop:2} leads to $\xi = 0$. 
\end{prf}

\textbf{Acknowledgment.} 
This note was written during a visit of the first-named author at
the Department of Mathematics of TU Darmstadt. 
The financial support and the excellent working conditions 
provided there are gratefully acknowledged.

\end{document}